\documentclass[12pt]{article}
\usepackage{amssymb}
\usepackage{latexsym}
\usepackage{amsfonts}

\def\LP{\mathrm{LP}}

\def\Rea{{\mathrm{Re~}}}
\def\Ima{{\mathrm{Im~}}}

\def\C{{\bf C}}    %added

\def\R{{\bf R}}    %added

\def\text{\hbox}   %added

\def\Im{\text{\rm Im~}}
\def\Re{\text{\rm Re~}}
\def\:{\colon}
\def\ds{\displaystyle\strut}

\title{Proof of a conjecture of
P\'olya on the zeros of successive derivatives of real entire functions}
\author{
Walter Bergweiler\thanks{Supported by the Alexander von
Humboldt Foundation, by the NSF, and by the G.I.F.,
the German--Israeli Foundation for Scientific Research and
Development, Grant G -809-234.6/2003.}~
and Alex Eremenko\thanks{
Supported by the NSF grants DMS-0555279 and DMS-0244547.}}
\begin{document}
\maketitle
\begin{abstract}
We prove P\'olya's conjecture of 1943:
For a real entire function of order greater than $2$
with finitely many non-real zeros, the number of non-real
zeros of the $n$-th derivative tends to infinity as $n\to\infty$.
We use the saddle point method and potential theory,
combined with the theory
of analytic functions with positive imaginary part
in the upper half-plane.
\end{abstract}

\noindent
{\bf 1. Introduction}
\vspace{.1in}

The theory of entire functions begins as a field of research in the
work of Laguerre \cite{Laguerre}, soon after the Weierstrass product
representation became available. Laguerre introduced the first
important classification
of entire functions, according to their genera.
We recall this notion. Let $f$ be an entire function and $(z_k)$ the sequence
of its zeros in $\C\backslash\{0\}$ repeated according to their multiplicities.
If $s$ is the smallest integer such that the series
$$\sum_{k=1}^\infty |z_k|^{-s-1}$$
converges, then $f$ has the Weierstrass representation
$$\displaystyle f(z)=z^me^{P(z)}\prod_{k=1}^\infty \left(1-\frac{z}{z_k}\right)
\exp\left(\frac{z}{z_k}+\ldots+\frac{1}{s}\left(\frac{z}{z_k}\right)^s
\right).$$
If $P$ in this representation
is a polynomial, then the genus is defined as $g:=\max\{ s,\deg P\}$.
If $P$ is a transcendental entire function,
or if the integer $s$ does not exist, then $f$ is said to have
infinite genus. A finer classification of entire functions by their
order
$$\rho=\limsup_{r\to\infty}\frac{\log\log M(r,f)}{\log r},\quad
\mbox{\rm where}\quad M(r,f)= \max_{|z|\leq r}|f(z)|,$$
was later introduced by Borel, based on the work of Poincar\'e and
Hadamard.
For functions of non-integral order
we have $g=[\rho]$,
and if $\rho$ is a positive integer, then $g=\rho-1$ or $g=\rho$.
So the order and the genus are simultaneously finite or infinite.

The state of the theory of entire functions in 1900 is described in
the survey of Borel \cite{Bo}. One of the first main problems of the theory was
finding relations between the zeros of a real entire function and the zeros of
its derivatives.
(An entire function is called real if it maps the real axis into
itself.)
If $f$ is a real polynomial with all zeros real, then all
derivatives $f^{(n)}$ have the same property. The proof based on Rolle's theorem
uses the finiteness of the set of zeros of $f$ in an essential way.
Laguerre discovered that this result still holds for entire functions
of genus $0$ or $1$, but in general fails
for functions of genus $2$ or higher \cite{Laguerre}.
P\'olya \cite{P1} refined this result as follows.
Consider the class of all entire functions which can be approximated
uniformly on compact subsets of the plane by real polynomials with all zeros real.
P\'olya proved that this class coincides with the set of functions
of the form 
$$f(z)=e^{az^2+bz+c}w(z),$$
where $a\leq 0$, $b$ and $c$ are real, and $w$ is a canonical product of
genus at most one, with all zeros real.
This class of functions is called the {\em Laguerre--P\'olya class}
and denoted by $\LP$. Evidently $\LP$ is closed with respect to differentiation,
so all derivatives of a function in $\LP$ have only real zeros. For a function
of class $\LP$, the order and the genus do not exceed $2$.

The class $\LP$ has several  other interesting characterizations, and it
plays an important role in many parts of analysis \cite{HW,K}.

The results of Laguerre and P\'olya on zeros of successive
derivatives inspired much research
in the 20th century.
In his survey article \cite{P4}, 
P\'olya writes: {\em ``The real axis
seems to exert an influence on the non-real
zeros of $f^{(n)}$;
it seems to attract these zeros when the order is less than
$2$, and it seems to repel them when
the order is greater than $2$''.} In the original text P\'olya wrote
``complex'' instead of ``non-real''. We replaced this in accordance
with the modern usage to avoid confusion. 
P\'olya then put this in a precise form by making the following two conjectures.
\vspace{.1in}

\noindent
{\bf  Conjecture A.} {\em
If the order of the real
entire function $f$ is less than
$2$, and $f$ has only a finite number
of non-real zeros, then its derivatives, from a certain one
onwards, will have no non-real zeros at all.
}
\vspace{.1in}

\noindent
{\bf  Conjecture B.} {\em If the order of the real
entire function $f$ is greater than
$2$, and $f$ has only a finite number
of non-real zeros, then the number of
non-real zeros of $f^{(n)}$ tends to infinity
as $n\to\infty$.}
\vspace{.1in}

Conjecture A was proved
by Craven, Csordas and Smith~\cite{CCS}. Later they proved in \cite{CCS1}
that the conclusion
holds if $\log M(r,f)=o(r^2)$. 
The result was refined by Ki and Kim~\cite{K1,KK}
who proved that the conclusion holds if $f=Ph$ where
$P$ is a real polynomial and $h\in\LP$. 

In this paper we establish Conjecture B.
Actually we prove that its conclusion holds if $f=Ph$, with a real
polynomial $P$ and a real entire function $h$ with real zeros that
does not belong to $\LP$.
\vspace{.1in}

\noindent
{\bf Theorem 1.} {\em Let $f$ be a real entire function of
finite order with
finitely many non-real zeros. Suppose that $f$ is not a product of
a real polynomial and a function of the class $\LP$.
Then the number $N(f^{(n)})$ of non-real zeros of $f^{(n)}$ satisfies
\begin{equation}
\label{card1}
\liminf_{n\to\infty}\frac{N(f^{(n)})}{n}>0.
\end{equation}
}
\vspace{.1in}

For real entire functions of infinite order with finitely many non-real zeros
we have the recent result
of Langley \cite{L} that $N(f^{(n)})=\infty$ for $n\geq 2$.
Theorem~1 and Langley's result together prove Conjecture B.
Combining Theorem~1 with the results
of Langley and Ki and Kim we conclude the following:
\vspace{.1in}

{\em For every real entire function $f$,
either $f^{(n)}$ has only real zeros for all sufficiently large $n$, or the number
of non-real zeros of $f^{(n)}$  tends to infinity with $n$.}
\vspace{.1in}

We give a short survey of the previous
results concerning Conjecture B.
These results can be divided into two
groups: asymptotic results on the zeros of $f^{(n)}$ as $n\to\infty$,
and non-asymptotic results for fixed~$n$.

The asymptotic results begin
with the papers of
P\'olya \cite{P2,P3}.
In the first paper, P\'olya proved that
if
\begin{equation}
\label{polya}
f=Se^T,
\end{equation}
where $S$ and $T$ are real polynomials, and $f^{(n)}$ has only real
zeros for all $n\geq 0$, then $f\in \LP$. In the second paper, he
introduced the {\em final set}, that is the
set of limit points of zeros of successive derivatives, and he found
this limit set for all meromorphic functions with at least two poles,
as well as for entire functions of the form (\ref{polya}).
It turns out that for a function
of the form (\ref{polya}),
the final set consists of equally
spaced rays in the complex plane. Namely,
if
$T(z)=az^d+bz^{d-1}+\ldots$, then
these rays emanate from the point $-b/(da)$ and have
directions $\arg z=(\arg a+(2k+1)\pi)/d,$ with $k=0,\ldots,d-1$.
It follows that the number of non-real zeros of $f^{(n)}$ for such a
function $f$ tends
to infinity as $n\to\infty$, unless $d=2$ and $a<0$, or $d\leq 1$,
that is, unless $\exp(T)\in\LP$. Notice that the final
set of a function (\ref{polya}) is independent of the polynomial $S$.

This result of P\'olya was generalized by McLeod \cite{Mc} to the case
that $S$ in (\ref{polya}) is an entire function satisfying
$\log M(r,S)=o(r^{d-1}),\; r\to\infty.$  This growth restriction seems natural:
a function $S$ of faster growth will influence the final set. This is seen
from P\'olya's result where the final set depends on~$b$.

All papers on final sets use some form of the saddle point method
to obtain an asymptotic expression
for $f^{(n)}$ when $n$ is large. This leads to conclusions
about the zeros of $f^{(n)}$.
McLeod used a
very general and powerful version of the saddle point method
which is due to Hayman \cite{H}. We do not survey here many
other interesting results on the final sets of entire functions,
as these results
have no direct bearing on Conjecture B, and we refer
the interested reader to \cite{B,Ge,Shen2} and references therein.

Passing to the non-asymptotic results, we need the following definition.
For a non-negative integer $p$ we define $V_{2p}$ as the class
of entire functions of the form
$$e^{az^{2p+2}}w(z),$$
where $a\leq 0$, and $w$ is a real entire function of genus at most $2p+1$
with all zeros real.
Then we define $U_0:=V_0$ and $U_{2p}:=V_{2p}\backslash V_{2p-2}$. 
Thus $\LP=U_0$.

Many of the non-asymptotic results were motivated by
an old conjecture, attributed to Wiman (1911) by
his student {\AA}lander \cite{Alan1,Alan}, that every
real entire function
$f$ such that $ff^{\prime\prime}$ has only real zeros belongs to the
class~$\LP$. For functions of finite order, Wiman made the more precise
conjecture that $f^{\prime\prime}$ has at least $2p$
non-real zeros if $f\in U_{2p}$.

Levin and Ostrovskii \cite{LO} proved that if $f$ is a real entire
function which satisfies
\begin{equation} \label{fast}
\limsup_{r\to\infty}\frac{\log\log M(r,f)}{r\log r}=\infty,
\end{equation}
then $ff''$  has infinitely many non-real zeros.
Hellerstein and Yang \cite{HY}
proved that under the same assumptions $ff^{(n)}$ has infinitely
many non-real zeros if $n\geq 2$. Thus Conjecture B was
established for functions satisfying (\ref{fast}).

A major step towards Conjecture B
was then made by Hellerstein and
Wil\-liam\-son \cite{HW1,HW2}. They proved that  for a real entire
function $f$, the condition that $ff'f^{\prime\prime}$ has only real
zeros implies that $f\in\LP$. It follows that for a real entire
function $f$ with only real zeros, such that $f\notin\LP$, the
number $N(f^{(n)})$ of non-real zeros of $f^{(n)}$ satisfies
$$\limsup_{n\to\infty}N(f^{(n)})>0.$$
The next breakthrough was made by Sheil-Small \cite{SS}
by proving Wiman's conjecture
for functions of finite order:
$f^{\prime\prime}$ has at least $2p$ non-real zeros if $f\in U_{2p}$. 
This was a new and deep result even for the case $f=e^T$ with
a polynomial $T$.
By refining the arguments of Sheil-Small, Edwards and Hellerstein
\cite{EH} were able to prove that
$$N(f^{(n)})\geq 2p,\quad\mbox{for}\quad n\geq 2\quad\mbox{and}\quad f=Ph,$$
where $P$ is a real polynomial and $h\in U_{2p}$. 

In the paper \cite{BEL}, the result of Sheil-Small
was extended to functions of infinite order, thus establishing
Wiman's conjecture in full generality. As we already mentioned,
Langley \cite{L} extended the result of \cite{BEL} to higher derivatives.
He proved the following:
Let $f$ be a real entire function of infinite
order with finitely many non-real zeros. Then  $N(f^{(n)})=\infty$
for all $n\geq 2$. Thus Conjecture B was established for functions
of infinite order.

To summarize the previous results, we can say that all asymptotic results 
were proved under strong a priori assumptions on the
asymptotic behavior of $f$, while the non-asymptotic results
estimated $N(f^{(n)})$ for $f\in U_{2p}$ from below in terms of $p$ rather than $n$.

In order to  state a refined version of our result,
we denote by $N_{\gamma,\delta}(f^{(n)})$ the number
of zeros of $f^{(n)}$ in 
$$\displaystyle\{ z: |\Ima z|> \gamma |z|, |z|>n^{1/\rho - \delta}\},$$
where $\rho$ is the order of $f$, and
$\gamma$ and $\delta$ are arbitrary positive numbers.
\vspace{.1in}

\noindent
{\bf Theorem 2}. {\em Suppose that $f=Ph$ where $P$ is a real polynomial
and $h\in U_{2p}$ with $p\geq 1$.
Then there exist
positive numbers
$\alpha$ and $\gamma$ depending only on $p$
such that
\begin{equation}
\label{card2}
\liminf_{n\to\infty}\frac{N_{\gamma,\delta}(f^{(n)})}{n}\geq \alpha>0,
\end{equation}
for every $\delta>0$.}
\vspace{.1in}

Estimates (\ref{card1}) and (\ref{card2}) seem to be new even for
functions of the form (\ref{polya}) with polynomials $S$ and $T$.
We recall that for 
a real entire function $f$ of genus $g$ with only finitely
many non-real zeros we have $N(f')\leq N(f)+g$
by a theorem of Laguerre and Borel \cite{Bo,G}.
Thus (\ref{card1}) and (\ref{card2}) give the right order of magnitude.
The question remains
how $\alpha$ and $\gamma$ in (\ref{card2}) depend  on~$p$.
\vspace{.1in}

We conclude this Introduction with a sketch of the proof of Theorem~1,
which combines the saddle point method used in the asymptotic 
results with potential theory and
the theory of analytic functions with positive
imaginary part in the upper half-plane. The latter theory
was a common tool in all non-asymptotic results since the
discovery of Levin and Ostrovskii \cite{LO} that
for a real entire function $f$ with real zeros, the logarithmic
derivative is a product of a real entire function
and a function which maps the upper
half-plane into itself (Lemma~4 below).

Our proof consists of three steps.
\vspace{.1in}

1. {\em Rescaling.} We assume for simplicity that all
zeros of $f$ are real.
The Levin--Ostrovskii representation 
gives $L:=f'/f=P_0\psi_0$, where
$P_0$ is a real polynomial of degree at least $2$,
and $\psi_0$ has non-negative imaginary
part in the upper half-plane
(Lemma~4). So we have a good control of the
behavior of $f'/f$ in the upper half-plane (Lemma~1).
For any sequence $\sigma$
of positive integers we can find a subsequence $\sigma'$ and
positive numbers $a_k$ such that for $k\to\infty,\; k\in\sigma'$,
the following limits exist in the upper half-plane (Section~3):
$$q(z)=\lim_{k\to\infty} a_k L(a_kz)/k,\quad\mbox{and}\quad
u(z)=\lim_{k\to\infty} (\log|f^{(k)}(a_kz)|-c_k)/k,$$
with an appropriate choice of real constants $c_k$. The second limit
makes sense in $L^1_{\mathrm{loc}}$, and $u$ is a subharmonic function
in the upper half-plane (Lemma~5).
If the condition (\ref{card1}) is not satisfied, that is the $f^{(k)}$ have few
zeros in the upper half-plane, then $u$ will be a harmonic function.
Our goal is to show that this is impossible.
\vspace{.1in}

2. {\em Application of the saddle point method to find a functional
equation for $u$.}
We express the $f^{(k)}$ as Cauchy integrals
over appropriately chosen circles
in the upper half-plane, and apply the saddle
point method to find asymptotics of these integrals
as $k\to\infty$ (Section~4).
The Levin--Ostrovskii representation and
properties of analytic functions with positive imaginary part
in the the upper half-plane
give enough information for the estimates needed
in the saddle point method.
As a result, we obtain an expression of $u$
in terms of $q$ in a {\em Stolz angle}
at infinity, by which 
we mean a region of the form
$\{ z:|z|>R,\, \varepsilon<\arg z<\pi-\varepsilon\}$ with $R>0$ and
$\varepsilon>0$.
The expression we obtain for $u$ is
\begin{equation}
\label{1}
u(z-1/q(z))=\Rea \int_i^z q(\tau)d\tau+\log |q(z)|.
\end{equation}
This equation, which plays
a fundamental role in our proof, can be derived
heuristically as follows:
When applying Cauchy's formula to obtain an expression for $f^{(k)}(a_kw)$,
it seems reasonable to write the radius of the circle in the form $a_kr_k$, so
that
\begin{eqnarray*}
f^{(k)}(a_kw)
&=&
\frac{k!}{2\pi i}
\int\limits_{|\xi|=a_kr_k}\frac{f(a_kw+\xi)}{\xi^k}\frac{d\xi}{\xi}\\ \\
&=&
\frac{k!}{2\pi i}\int\limits_{|\zeta|=r_k}
\exp\left( \log f(a_k(w+\zeta))-k\log a_k\zeta\right)
\frac{d\zeta}{\zeta}.
\end{eqnarray*}
The saddle point method of finding
the asymptotic behavior of such an integral as $k\to\infty$ involves
a stationary point
of the function in the exponent,
that is
a solution of the equation
$$0=\frac{d}{d\zeta}\left(\frac{\log f(a_k(w+\zeta))}{k}-\log a_k\zeta\right)
=
\frac{a_kL(a_k(w+\zeta))}{k}-\frac{1}{\zeta}
=: q_k(w+\zeta)-\frac{1}{\zeta}.$$
This suggests to take $r_k=|\zeta|$ where $\zeta$ is a solution
of the last equation.
Setting $z=w+\zeta,$ we obtain $w=z-1/q_k(z),$
and the saddle point method gives
$$\frac{1}{k}\left(\log f^{(k)}\left(a_k 
(z-1/q_k(z))\right)-c_k\right)\sim
\int_i^z q_k(\tau)d\tau+\log q_k(z),$$
with some constants $c_k$. From this we derive (\ref{1}).
\vspace{.1in}

3. {\em Study of the functional equation $(\ref{1})$.}
If the $f^{(k)}$ have few zeros in the upper half-plane, the function
$u$ should be harmonic in the upper half-plane. On the other hand,
we show that a harmonic function $u$ satisfying (\ref{1})
in a Stolz angle
cannot have a harmonic continuation into the whole upper half-plane
(Section~5). Here again we use the properties of analytic functions
with positive imaginary part in the upper half-plane.
\vspace{.1in}

The authors thank David Drasin, Victor Katsnelson, Jim Langley,
Iosif Ostrovskii and the referees for useful comments on this work.
\vspace{.2in}

\noindent
{\bf 2. Preliminaries}
\vspace{.2in}

In this section we collect for the reader's convenience
all necessary facts on potential theory and on functions
with positive imaginary part in the upper half-plane 
$H:=\{ z:\Ima z>0\}.$
\vspace{.1in}

\noindent
{\bf Lemma 1.} {\em Let $\psi\not\equiv 0$ be an analytic function in 
$H$
with non-negative imaginary part.
Then
\begin{equation}
\label{21}
|\psi(i)|\frac{\Im z}{(1+|z|)^2}
\leq|\psi(z)|\leq |\psi(i)| \frac{(1+|z|)^2}{\Im z},
\end{equation}

\begin{equation}
\label{22}
\left|\frac{\psi'(z)}{\psi(z)}\right|\leq\frac{1}{\Ima z},
\end{equation}
and
\begin{equation}
\label{23}
\left|\log\frac{\psi(z+\zeta)}{\psi(z)}\right|\leq 1
\quad\mbox{for}\quad |\zeta|<\frac{1}{2}\Ima z.
\end{equation}
}
\vspace{.1in}

Inequality (\ref{21}) is a well-known estimate due to
Carath\'eodory, see for example \cite[\S 26]{V}. 
Inequality (\ref{22})
is the Schwarz Lemma: it says that the derivative of $\psi$
with respect to the hyperbolic metric in $H$ is at most $1$.
Finally (\ref{23}) follows from
(\ref{22}) by integration.
\vspace{.1in}

\noindent
{\bf Lemma 2.} {\em A holomorphic function $\psi$ in $H$
has non-negative imaginary part 
if and only if it has the form
$$\psi(z)=a+\lambda z+\int_{-\infty}^{\infty}
\left(\frac{1}{t-z}-\frac{t}{1+t^2}\right)d\nu(t),$$
where $\lambda\geq 0$, $a$ is real, and
$\nu$ is a non-decreasing
function of finite variation on the real line.}
\vspace{.1in}

In the case that $\psi$ is a meromorphic function
in the whole plane, $\nu$
is piecewise constant with jumps at the poles of $\psi$.
Then the representation in Lemma~2 is similar to
the familiar Mittag--Leffler representation.
Lemma~2 can be found in \cite{KacK}, where
it is derived from the similar Riesz--Herglotz representation
of functions with positive imaginary part in the unit disc.
\vspace{.1in}

\noindent
{\bf Lemma 3.} {\em Let $G$ be a
holomorphic function in $H$ with non-negative
imaginary part. Then there exists $\lambda\geq 0$ such that
$$G(z)=\lambda z+G_1(z),$$
where $\Im G_1(z)\geq 0,\; z\in H$, and $G_1(z)=o(z)$
as $z\to\infty$ in every Stolz angle.}
\vspace{.1in}

This is a theorem of Wolff; see,
for example, \cite[\S 26]{V},
where it is derived from the Schwarz Lemma.
Another proof can be easily
obtained from Lemma~2.
The number $\lambda$ is called the
{\em angular derivative} of $G$ at infinity.
\vspace{.1in}

\noindent
{\bf Lemma 4.} {\em Let $h$ be a function of class $U_{2p}$.
Then the logarithmic derivative of $h$ has
a representation
\begin{equation}
\label{lo}
h'/h=P_0\psi_0,
\end{equation}
where $P_0$ is a real polynomial, $\deg P_0= 2p$,
the leading coefficient of $P_0$ is
negative, and $\psi_0\not\equiv 0$ is a function with non-negative
imaginary part in $H$.}
\vspace{.1in}

Equation (\ref{lo}) is the Levin--Ostrovskii
representation already mentioned, which was used in many papers
on the subject \cite{BEL,EH,HW1,HY,L,LO,SS}.
Levin and Ostrovskii \cite{LO}
did not assume that $h$ has finite order
and showed that (\ref{lo}) holds with some real entire
function $P_0$. Hellerstein and Williamson \cite[Lemma 2]{HW1}
showed that $P_0$
is a polynomial if $h$ has finite order, and they
gave estimates of the degree of $P_0$ in \cite[Lemma 8]{HW1}.
As they noted, an upper bound for the degree of $P_0$ 
follows from an old result of Laguerre \cite[p. 172]{Laguerre}.

The factorization (\ref{lo}) is not unique, and the degree of $P_0$
is not uniquely determined by $h$. As our version of this
factorization is different from
\cite[Lemma 8]{HW1}, we include
a proof. 
\vspace{.1in}

{\em Proof of Lemma 4.} 
First we consider the simple case that $h=e^T$.

If $\deg T=2p+2$, then the leading coefficient
of $T$ is negative by definition of the class $U_{2p}$.
As the degree of $T'$ is odd, $T'$ has a real root $c$.
So we can put $P_0(z):=T'(z)/(z-c)$ and $\psi_0(z):=(z-c)$.

If $\deg T=2p+1$, then we set $P_0:=\pm T'$ and $\psi_0:=\pm 1$,
where the sign is chosen to ensure that the leading coefficient
of $T'$ is negative.

Finally, if $\deg T=2p$, then the condition $h\in U_{2p}$
implies that the leading coefficient of $T$ is positive,
and we set $P_0(z):=-zT'(z)$ and $\psi_0(z):=-1/z$.

From now on we assume that $h$ has at least one real
zero $a_0$.
If $h$ has only finitely many negative zeros
and $h(x)\to 0$ as $x\to
-\infty$, then
we also consider $-\infty$ as a zero of $h$. Similarly we consider
$+\infty$ as a zero of $h$
if $h$ has only finitely many positive zeros and
$h(x)\to 0$ as $x\to +\infty$.
We arrange the zeros of $h$ into an increasing sequence $(a_j)$,
where each zero occurs once, disregarding multiplicity. 
The range of the subscript $j$ will be
$M<j<N$, where
$-\infty\leq M< 0\leq N\leq \infty$,
with $a_{M+1}=-\infty$ and $a_{N-1}=+\infty$ in the cases described
above.

By Rolle's theorem, each open interval $(a_j,a_{j+1})$
contains a zero $b_j$
of $h'$. To make a definite choice,
we take for $b_j$ the
largest or the smallest zero in this interval.
Each $b_j$ occurs in this sequence
only once, and we disregard multiplicity. 
We define
$$\psi(z):=\frac{1}{z-a_{N-1}}\prod_{M<j<N-1}\frac{1-z/b_j}{1-z/a_j}$$
where the
factor $z-a_{N-1}$ is omitted
if $a_{N-1}=+\infty$ or $N=\infty$,
and the factor $1-z/a_{M+1}$ is omitted if
$a_{M+1}=-\infty$. If for some $j\in (M,N-1)$ we have $a_j=0$
or $b_j=0$ then the $j$-th  factor has to be replaced
by $(z-b_j)/(z-a_j)$. As in \cite{LO,Levin}
we see that the product converges and
is real only on the real axis.
We define $P:=h'/(h\psi)$.
Then standard estimates using the Lemma on the
Logarithmic Derivative and Lemma~1 imply that $P$ is a polynomial
\cite{HW1,LO}.
In particular, $P$ has only finitely many zeros.
The zeros of $P$ are precisely the zeros
that $h'$ has in addition
to the $b_j$ and possible multiple zeros $a_j$ of $h$.
These additional zeros were called {\em extraordinary}
zeros by {\AA}lander \cite{Alan1}.
An extraordinary zero may be real or not,
it can be real but different from all $b_j$, or may be one
of the $b_j$: if $b_j$ is a zero of $h'$ of multiplicity $n\geq 2$,
then $b_j$ is considered an extraordinary zero of multiplicity $n-1$.
Since the number of zeros of $h'$ in every interval $(a_j,a_{j+1})$ is
odd (counted with multiplicity), the number of extraordinary zeros
of $h'$ in this interval is even. As the non-real extraordinary zeros come in
complex conjugate pairs, their number is also even.
Overall, $h'$ has an even number of extraordinary zeros,
counted with multiplicity.  Thus the degree of $P$ is even.

We choose $\varepsilon=\pm 1$ so that $P_0:=\varepsilon P$ has 
negative leading coefficient and define $\psi_0:=\varepsilon\psi$.
Then (\ref{lo}) holds.
Since the number of extraordinary zeros in each interval
$(a_j,a_{j+1})$ is even, $P_0$ has an even number
of zeros to the left of $a_0$, and, as $P_0(x)\to-\infty$
as $x\to-\infty$, we conclude that 
$P_0(a_0)<0$.
Since $h'/h$ has a simple pole with
positive residue at $a_0$, and $P_0(a_0)<0$,
we obtain that $\psi_0$ has a pole with negative residue at $a_0$.
Since $\psi_0$ is real only on the real axis this implies that
$\psi_0(H)\subset H$.
We note that
$$\psi_1:=-1/\psi_0$$ also maps $H$ into itself. 

Next we show that $\deg P_0\leq 2p$.
We recall that $h$ has the
form 
\begin{equation}
\label{xxx}
h(z)=e^{az^{2p+2}}w(z),
\end{equation}
where $a\leq 0$, and $w$ is a real entire
function of genus at most $2p+1$.
Then
$\log M(r,w)=o(r^{2p+2})$ as $r\to\infty$; see, for example,
\cite[I.4]{Levin}. 
Schwarz's 
formula \cite[I.6]{Levin}
shows that $|w'(z)/w(z)|=o(z^{2p+1})$ as $z\to\infty$ in every
Stolz angle. We thus have $h'(z)/h(z)= b z^{2p+1}+o(z^{2p+1})$
as $z\to\infty$ in every
Stolz angle, with $b:=(2p+2)a\leq 0$.
Denote by $\lambda_1$ the angular derivative of $\psi_1$ (see Lemma~3).
Then $\psi_1(z)=\lambda_1 z +o(z)$ as $z\to\infty$ in every
Stolz angle. Altogether we find that 
$$P_0(z)=-\frac{h'(z)}{h(z)}\psi_1(z)=
-b \lambda_1 z^{2p+2}+o(z^{2p+2}),
$$
as $z\to\infty$ in every Stolz angle. 
If $b<0$ and $\lambda_1>0$, then
$P_0$ has the positive leading coefficient 
$-b\lambda_1$, a contradiction.
Thus $b=0$ or $\lambda_1=0$ so that $|P_0(z)|=o(z^{2p+2})$
and hence $\deg P_0< 2p+2$. Since $\deg P_0$ is even we
thus obtain  $\deg P_0\leq 2p$.

Now we show that $\deg P_0\geq 2p$.
We have 
$$P_0(z)=cz^d+\ldots,\quad\mbox{where}\quad c<0.$$
Let $\lambda_0\geq 0$ be the angular derivative of
$\psi_0$. Then in every Stolz angle
$$\frac{h'(z)}{h(z)}=c\lambda_0z^{d+1}+o(z^{d+1}),\quad z\to\infty.$$
Integrating this along straight lines, we conclude
that 
$$\log h(z)=\frac{c\lambda_0}{d+2}z^{d+2}+o(z^{d+2}),\quad z\to\infty.$$
If $c\lambda_0<0$, we compare this with (\ref{xxx}) 
and obtain that $d=2p$.
If $\lambda_0=0$ we obtain that $a=0$ in (\ref{xxx}),
so the genus $g$ of $h$ is at most $2p+1$.

Now we follow \cite[Lemma 8]{HW1}.
The logarithmic derivative of $h$ has the form
\begin{equation}
\label{raz}
\frac{h'(z)}{h(z)}=Q(z)+z^g\sum_j\frac{m_j}{a_j^g(z-a_j)},
\end{equation}
where $g\leq 2p+1$ is the genus of $h$,
$m_j$ the multiplicity of the zero $a_j$
and $Q$ a polynomial. On the other hand,
(\ref{lo}) combined with Lemma~2 gives
\begin{equation}
\label{dwa}
\frac{h'(z)}{h(z)}=P_0(z)\left\{ \lambda_0z+c_0+\sum_j A_j\left(\frac{1}{a_j-z}-
\frac{a_j}{1+a_j^2}\right)\right\},
\end{equation}
where $\lambda_0\geq 0$, $A_j\geq 0$ and $c_0$ is real. We also
have
\begin{equation}
\label{3.54}
\sum_j \frac{A_j}{a_j^2}<\infty.
\end{equation}
Equating the residues at each pole in the expressions
(\ref{raz}) and (\ref{dwa}) 
we obtain
\begin{equation}
\label{3.57}
P_0(a_j)=-m_j/A_j<0.
\end{equation}
Now we choose $C>0$ such that
$0<-P_0(a_j)\leq C|a_j|^d$. Then (\ref{3.54}) and (\ref{3.57})
imply
$$\frac{1}{C}\sum_j\frac{m_j}{|a_j|^{d+2}}\leq\sum_j
\frac{-m_j}{a_j^2P_0(a_j)}
=\sum_j\frac{A_j}{a_j^2}<\infty,$$
which shows $d+2\geq g+1$, that is $d\geq g-1$. If $g=2p+1$, we have
$d\geq 2p$. If $g=2p$, then also $d\geq 2p$ because $d$ is even.
It remains to notice that one cannot
have $g\leq 2p-1$, because $h$ is of genus $g$ and belongs to
$U_{2p}$. This completes the proof of Lemma~4. 

\vspace{.1in}

\noindent
{\bf Lemma 5.} {\em Let $(u_k)$ be a sequence of subharmonic functions
in a region $D$, and suppose that the $u_k$ are uniformly bounded from above
on every compact subset of $D$. Then one can choose a subsequence
of $(u_k)$, which either converges to $-\infty$
uniformly on compact subsets of $D$ or
converges in $L^1_{\mathrm{loc}}$ (with respect to
the Lebesgue measure in the plane)
to a subharmonic
function $u$. In the latter case,
the Riesz measures of the $u_k$
converge weakly to the Riesz measure of $u$.}
\vspace{.1in}

This result can be found in 
\cite[Theorem 4.1.9]{Hor}. 
We recall that the Riesz measure of a subharmonic function $u$
is $(2\pi)^{-1}\Delta u$ in the sense of distributions.
\vspace{.2in}

{\bf 3. Beginning of the proof: rescaling}
\vspace{.2in}

We write
$$f=Ph,$$
where $P$
is a real polynomial and $h$ is a real entire function of finite order, with
all zeros real, and $h\notin\LP$. 
By Lemma 4, we have 
$$h'/h=P_0\psi_0,$$
where $\psi_0$ is as in  Lemma~4,
and $P_0$ is a polynomial of degree $2p$ 
whose leading coefficient is negative.
Note that $p\geq 1$ since $h\notin \LP$. We have
\begin{equation}
\label{3a}
f'/f=P_0\psi_0+P'/P.
\end{equation}
Using (\ref{21}) we obtain
\begin{equation}
\label{3b}
\frac{rf'(ir)}{f(ir)}\to\infty\quad\mbox{as}\quad r\to\infty,\; r>0.
\end{equation}

Fix $\gamma,\delta>0$, and let $\sigma$ be
a sequence of positive integers
along which the lower limit in (\ref{card2}) is attained; that is,
$$\beta:=\liminf_{k\to\infty}\frac{N_{\gamma,\delta}(f^{(k)})}{k}=
\lim_{k\to\infty,\; k\in\sigma}\frac{N_{\gamma,\delta}(f^{(k)})}{k}.$$
In the course of the proof we will choose subsequences
of $\sigma$, and will continue to denote them by the same letter $\sigma$.

By (\ref{3b}), for every large 
$k\in \sigma$, there exist $a_k>0,\; a_k\to\infty,$ such that
$$\left|\frac{a_kf'(ia_k)}{f(ia_k)}\right|=k.$$
An estimate of the logarithmic derivative using Schwarz's formula implies that
$$\left|\frac{a_kf'(ia_k)}{f(ia_k)}\right|\leq \left|a_k\right|^{\rho+o(1)}\quad
\mbox{as}\quad k\to\infty,$$ where $\rho$ is the order of $f$.
We may thus assume that
\begin{equation}
\label{a_k}
|a_k|\geq k^{1/\rho-\delta/2}
\end{equation}
for $k\in\sigma$.
We define
\begin{equation}
\label{qk}
q_k(z):=\frac{a_kf'(a_kz)}{kf(a_kz)}.
\end{equation}
Then $|q_k(i)|=1$. From (\ref{3a}) and (\ref{21})
we deduce that the $q_k$  are uniformly bounded on compact subsets of $H$,
and thus the $q_k$ form a normal family in $H$. Passing to a subsequence,
we may assume that
\begin{equation}
\label{3c}
q_k\to q,
\end{equation}
as $k\to\infty$, $k\in\sigma$,
uniformly on compact subsets in $H$.
We choose a branch of the logarithm in a neighborhood of $f(ia_k)$,
put $b_k:=\log f(ia_k)$, and define
\begin{equation}
\label{defQ}
Q(z):=\int_i^zq(\zeta)d\zeta
\end{equation}
and 
$$Q_k(z):=\int_i^zq_k(\zeta)d\zeta.$$
It follows from (\ref{3c}) that
$$Q_k(z)=\frac{1}{k}\left(\log f(a_kz)-b_k\right)\to Q(z),\quad
k\to\infty,\; k\in \sigma,$$
uniformly on compact subsets of $H$, and $Q_k(i)=0$.
The chosen branches of the $\log f$ are well defined on every
compact subset of $H$, if $k$ is large enough,
because $f$ has only finitely
many zeros in the upper half-plane, and $a_k\to\infty$.

Let $z$ be a point in $H$, and $0<t<\Ima z$. Then the disc
$\{\zeta:|\zeta-a_kz|<ta_k\}$ is contained in $H$ and does not contain
any zeros of $f$ if $k$ is large enough. Thus, by Cauchy's formula,
\begin{eqnarray}
\displaystyle
\nonumber
f^{(k)}(a_kz)&=&\frac{k!}{2\pi i}\int_{|\zeta|=a_kt}
\frac{f(a_kz+\zeta)}{\zeta^k}\frac{d\zeta}{\zeta}\nonumber\\
\label{3d}\\
&=&\frac{k!}{2\pi i}\int_{|\zeta|=a_kt}\frac{\exp\left(kQ_k(z+\zeta/a_k)+b_k\right)}{
\zeta^k}\frac{d\zeta}{\zeta}.
\nonumber
\end{eqnarray}
So
$$|f^{(k)}(a_kz)|\leq\frac{k!}{(a_kt)^k}\exp\left(\Rea b_k+k\max_{|\zeta|=t}\Rea Q_k(z+\zeta)\right),$$
and
with
\begin{equation}
u_k(z):=\frac{1}{k}\left(\log|f^{(k)}(a_kz)|-\Rea b_k-\log k!\right)+
\log a_k
\label{defuk}
\end{equation}
we obtain
\begin{equation}
\label{bounduk}
u_k(z)\leq\max_{|\zeta|=t}\Rea Q_k(z+\zeta)-\log t.
\end{equation}
Since $Q_k\to Q$, we deduce that the $u_k$ are uniformly bounded from above
on
compact subsets of $H$. By Lemma~5, after choosing a subsequence,
we obtain
\begin{equation}
\label{converg}
u_k\to u,
\end{equation}
where $u$ is a subharmonic function
in $H$ or $u\equiv-\infty$. The convergence in (\ref{converg})
holds in $H$ in the sense described in Lemma~5.
We will later see that $u\not\equiv-\infty$.
We will then show that $u$ cannot be harmonic in $H$. This will
prove Theorem~1.

Moreover, 
we will see that there exist positive constants $\gamma$ and
$\alpha$ depending only on $p$,
such that the total Riesz measure of $u$
in the region  $\{ z: \Ima z> \gamma |z|\}$ is 
at least $\alpha/2.$

On the other hand, it follows from the definition of
$u_k$ and $N_{\gamma,\delta}(f^{(k)})$
that $N_{\gamma,\delta}(f^{(k)})/(2k)$ is 
the total Riesz measure of $u_k$ in the region
$\{ z: \Ima z> \gamma |z|, |z|>k^{1/\rho - \delta}/a_k \}.$
Note that 
$k^{1/\rho - \delta}/a_k \to 0$
by  (\ref{a_k}).
Passing to the limit as $k\to\infty$, $k\in\sigma$, we obtain
$\beta\geq \alpha>0$,
which will complete the proof of Theorem~2.
\vspace{.2in}

\noindent
{\bf 4. Application of the saddle point method}
\nopagebreak
\vspace{.2in}

It follows from (\ref{3a}) and the definitions of $q_k$ in (\ref{qk})
and $q$ in (\ref{3c}) that
\begin{equation}
\label{4a}
q(z)=-z^{2p}\psi(z),
\end{equation}
where
$\psi:H\to\overline{H}\backslash\{0\}$ is defined by
$$
\psi(z):=\lim_{k\to\infty}\frac{\psi_0(a_kz)}{|\psi_0(i a_k)|}.
$$
Here $\psi_0$ is the real meromorphic function in $\C$ from (\ref{3a}).
Note that
\begin{equation}
\label{i}
|q(i)|=|\psi(i)|=1.
\end{equation}

Now we define
$$F(z):=z-1/q(z).$$
{}From (\ref{4a}) and (\ref{21}), it follows
that $1/q(z)\to 0$ as $z\to\infty$ in every Stolz angle.
So
$$F(z)\sim z,\quad z\to\infty$$
in every Stolz angle.
We fix $\delta_0\in (0,\pi/2)$, for example $\delta_0=\pi/4$.
It follows
that there exists $R>0$ such that
$\Ima F(z)>0$ in the region
\begin{equation}
\label{defSR}
S_{R}:=\{ z:|z|>R,\; \delta_0<\arg z<\pi -\delta_0\}.
\end{equation}
Moreover,
if $\delta_1\in(0,\delta_0)$, then
$$F(z)\in S':=\{ w:\delta_1<\arg z<\pi-\delta_1\},\quad\mbox{when}
\quad z\in S_{R},$$
provided that $R$ is large enough, say $R>R_0$.
We notice that $R_0$ is independent of $f$; this follows
from (\ref{21}), (\ref{4a}) and (\ref{i}). We will enlarge
$R_0$ in the course of the proof but it will be always a constant
independent of $f$. 
We will obtain for $R>R_0$ the identity
\begin{equation}
\label{identity}
u(F(z))=\Rea Q(z)+\log|q(z)|,\quad z\in S_R.
\end{equation}

In order to do this, we use Cauchy's formula (\ref{3d}) with $z$ replaced by
$F_k(z):=z-r_k(z)$, and $t=|r_k(z)|,$
where we have set
$r_k:=1/q_k$ to simplify our formulas.

We obtain
\begin{eqnarray}\ds
f^{(k)}(a_kF_k(z))&=&\frac{k!}{2\pi}\int_{-\pi}^\pi
\frac{\exp\left(kQ_k(F_k(z)+r_k(z)e^{i\theta})+b_k\right)}{
\left(a_kr_k(z)e^{i\theta}\right)^k}d\theta\nonumber\\
\label{4a1}\\
&=&\frac{k!q_k^k(z)}{2\pi a_k^k}e^{b_k}\int_{-\pi}^\pi
\exp\left(kQ_k(F_k(z)+r_k(z)e^{i\theta})-ik\theta\right)d\theta.
\nonumber
\end{eqnarray}

To determine the asymptotic behavior of this integral as $k\to\infty$,
we expand
$$L_k(\theta):=Q_k(F_k(z)+r_k(z)e^{i\theta})$$
into a Taylor series:
$$L_k(\theta)=L_k(0)+L^\prime_k(0)\theta+\frac{1}{2}L_k^{\prime\prime}(0)
\theta^2+\frac{1}{6}E_k(\theta)\theta^3,$$
where
$$|E_k(\theta)|\leq\max_{t\in[0,\theta]}\left|L_k^{\prime\prime\prime}(t)
\right|.$$
We notice that
$$L_k(0)=Q_k(F_k(z)+r_k(z))=Q_k(z),$$
and
$$
L_k^\prime(\theta)=iQ_k^\prime(F_k(z)+r_k(z)e^{i\theta})r_k(z)e^{i\theta}
= iq_k(F_k(z)+r_k(z)e^{i\theta})r_k(z)e^{i\theta},
$$
so that $L_k^\prime(0)=i$.
Moreover,
$$
L_k^{\prime\prime}(\theta)=
-q_k^\prime(F_k(z)+r_k(z)e^{i\theta})r_k^2(z)e^{2i\theta}
-q_k(F_k(z)+r_k(z)e^{i\theta})r_k(z)e^{i\theta},
$$
so that
\begin{equation}
\label{Lpp}
L_k^{\prime\prime}(0)=-\frac{q_k^\prime(z)}{q^2_k(z)}-1.
\end{equation}
Finally, we have
\begin{eqnarray*}
L_k^{\prime\prime\prime}(\theta)&=&-iq_k^{\prime\prime}(F_k(z)+r_k(z)e^{i\theta})
r_k^3(z)e^{3i\theta}\\ \\
&&-3iq_k^{\prime}(F_k(z)+r_k(z)e^{i\theta})r_k^2(z)e^{2i\theta}\\ \\
&&-iq_k(F_k(z)+r_k(z)e^{i\theta})r_k(z)e^{i\theta}.
\end{eqnarray*}
To estimate $L_k^{\prime\prime}(0)$ and $L_k^{\prime\prime\prime}(\theta)$
we notice that
$$\frac{q'(\zeta)}{q(\zeta)}=\frac{2p}{\zeta}+\frac{\psi^\prime(\zeta)}{\psi
(\zeta)}.$$
It follows from  (\ref{22}) that
$$
\left|\frac{q'(\zeta)}{q(\zeta)}\right|\leq\frac{2p+1}{\Ima\zeta}.
$$
{}From this we deduce that
\begin{eqnarray*}
\left|\frac{d}{d\zeta}
\left(\frac{q^{\prime}(\zeta)}{q(\zeta)}\right)\right|&=&
\frac{1}{2\pi}\left|\int_{|z-\zeta|=\frac12 \Ima\zeta}
\frac{q'(z)}{q(z)(z-\zeta)^2}dz \right|\\
&\leq& \frac{2}{\Ima\zeta}\max_{|z-\zeta|=\frac12 \Ima\zeta} 
\left|\frac{q'(z)}{q(z)}\right|\\
&\leq& 
\frac{4(2p+1)}{(\Ima\zeta)^2}.
\end{eqnarray*}
so that
$$
\left|\frac{q^{\prime\prime}(\zeta)}{q(\zeta)}\right|
=
\left|\frac{d}{d\zeta}
\left(\frac{q^{\prime}(\zeta)}{q(\zeta)}\right)+
\left(\frac{q^{\prime}(\zeta)}{q(\zeta)}\right)^2\right|
\leq
\frac{(2p+5)(2p+1)}{(\Ima\zeta)^2}.
$$
Since $q_k\to q$ 
we deduce from the above estimates that 
\begin{equation}
\label{4b2}
\left|\frac{q_k'(\zeta)}{q_k(\zeta)}\right|\leq\frac{2p+2}{\Ima\zeta}
\ \ \ \mbox{\rm and}\ \ \
\left|\frac{q_k^{\prime\prime}(\zeta)}{q_k(\zeta)}\right|\leq
\frac{(2p+5)^2}{(\Ima\zeta)^2}
\end{equation}
on any compact subset of $H$, provided $k$ is sufficiently large, 
$k\in\sigma$.

Fix $\eta>0$. It follows from (\ref{Lpp}) and (\ref{4b2}) 
that if $z\in S_{R}$, where $R>R_0$ and $k$ is large enough,
then
$$|L_k^{\prime\prime}(0)+1|<\eta.$$
In particular,
\begin{equation}
\label{4c}
\Rea L_k^{\prime\prime}(0)\leq -1+\eta.
\end{equation}
Moreover, if $z\in S_{R}$ with $R>R_0$,
then
\begin{equation}
\label{zeta}
\zeta:=F_k(z)+r_k(z)e^{i\theta}\in S',\quad\mbox{and}\quad\zeta\to\infty\quad
\mbox{as}\quad z\to\infty.
\end{equation}
For large $|z|$ we have
$$
|L_k^{\prime\prime\prime}(\theta)|\leq\left|\frac{q_k(\zeta)}{q_k(z)}
\right|\left(\frac{|q_k^{\prime\prime}(\zeta)|}{|q_k(\zeta)||q_k(z)|^2}
+
3\left|\frac{q_k^\prime(\zeta)}{q_k(\zeta)q_k(z)}\right|+1\right).
$$
By (\ref{23}) and (\ref{4a}) we have
$$q_k(\zeta)/q_k(z)\to1,$$
where $\zeta$ is defined in (\ref{zeta}),
uniformly with respect to $z\in S_{R}$ 
and $\theta\in[-\pi,\pi]$ as $R\to\infty$.
Combining these estimates with (\ref{4b2}) we find that
$|L_k^{\prime\prime\prime}(\theta)|\leq 1+\eta$ and hence
\begin{equation}
\label{4d}
|E_k(\theta)|\leq 1+\eta,
\end{equation}
for $z\in S_{R}$ and $|\theta|\leq \pi$, provided $R>R_0$.
Altogether we have the Taylor expansion
\begin{equation}
\label{4e}
L_k(\theta)=Q_k(z)+i\theta+\frac{1}{2}L_k^{\prime\prime}(0)\theta^2+
\frac{1}{6}E_k(\theta)\theta^3,
\end{equation}
with $L_k^{\prime\prime}(0)$ and $E_k(\theta)$ satisfying
(\ref{4c}) and (\ref{4d}).
We define $C_k(z):=-\frac12 L_k^{\prime\prime}(0)$ and notice that in view of
(\ref{Lpp})
$$C_k(z)\to C(z):=\frac{1}{2}\left(1+\frac{q^\prime(z)}{q^2(z)}\right),\quad k\to\infty,k\in\sigma,$$
uniformly on any compact subset in $H$.

We note that if
$$|\theta|\leq\theta_0:=\frac{3(1-3\eta)}{1+\eta},$$
then
\begin{eqnarray}
& & \Rea\left(-C_k(z)
\theta^2+\frac{1}{6}E_k(\theta) \theta^3\right)\nonumber\\
&=&
\Rea\left(\frac{1}{2}L_k^{\prime\prime}(0)
\theta^2+\frac{1}{6}E_k(\theta) \theta^3\right)\nonumber\\
&\leq&\theta^2\left(\frac{1}{2}(-1+\eta)+\frac{1}{6}(1+\eta)|\theta|\right)
\label{dominated}\\
&=&\theta^2\left(\frac{1}{2}(-1+\eta)+\frac{1}{2}(1-3\eta)\right)\nonumber\\
&\leq&-\eta\theta^2.\nonumber
\end{eqnarray}
Using (\ref{4e}) we obtain
\begin{eqnarray*}
\ds
&{}&\int_{-\theta_0}^{\theta_0}\exp\left(
kQ_k(F_k(z)+r_k(z)e^{i\theta})-ik\theta\right)d\theta\\
&=&e^{kQ_k(z)}\int_{-\theta_0}^{\theta_0}\exp\left(
k\left(-C_k(z)\theta^2+\frac{1}{6}E_k(\theta)\theta^3\right)\right)d\theta\\
&=&\frac{e^{kQ_k(z)}}{\sqrt{k}}\int_{-\theta_0\sqrt{k}}^{\theta_0\sqrt{k}}\exp\left(
-C_k(z)t^2+\frac{1}{6}E_k(t/\sqrt{k})\frac{t^3}{\sqrt{k}}\right)dt
\end{eqnarray*}
Combining this with
(\ref{dominated}) and the theorem on dominated convergence
we obtain
\begin{equation}
\label{4e1}
\int_{-\theta_0}^{\theta_0}\exp\left(
kQ_k(F_k(z)+r_k(z)e^{i\theta})-ik\theta\right)d\theta
\sim\frac{e^{kQ_k(z)}}{\sqrt{k}}\sqrt{\frac{\pi}{C(z)}}
\end{equation}
as $k\to\infty$, 
$k\in\sigma$,
uniformly on compact subsets of $S_{R}$, for $R>R_0$.

In order to estimate the rest of the integral (\ref{4a1}) we will show
that
\begin{equation}
\label{4s}
\Rea Q_k(F_k(z)+r_k(z)e^{i\theta})\leq\Rea Q_k(z)-\frac{1}{2}(1-\cos\theta_0),
\end{equation}
for $\theta_0\leq\theta\leq\pi$. We have
\begin{eqnarray*}\ds
&{}&Q_k(z)-Q_k(F_k(z)+r_k(z)e^{i\theta})\\ \\
&=&Q_k(F_k(z)+r_k(z))-Q_k(F_k(z)+r_k(z)e^{i\theta})\\ \\
&=&\int_{r_k(z)e^{i\theta}}^{r_k(z)} q_k(F_k(z)+\zeta)d\zeta\\ \\
&=&\int_0^1 q_k(F_k(z)+r_k(z)e^{i\theta}+tr_k(z)(1-e^{i\theta}))
(1-e^{i\theta})r_k(z)dt\\ \\
&=&(1-e^{i\theta})\int_0^1\frac{q_k(\zeta_t)}{q_k(z)}dt,
\end{eqnarray*}
where
$$\zeta_t:=F_k(z)+r_k(z)e^{i\theta}+tr_k(z)(1-e^{i\theta}).$$
Now $\zeta_t/z\to 1$ as $z\to\infty,\; z\in S_{R}$, uniformly
with respect to $t\in [0,1]$. Using (\ref{23}) we see that $q_k(\zeta_t)/q_k(z)\to 1$
for $z\in S_{R}$ as $R\to\infty.$ In particular, we have
$$\Rea\left((1-e^{i\theta})\frac{q_k(\zeta_t)}{q_k(z)}\right)\geq
\frac{1-\cos\theta}{2}\geq\frac{1-\cos\theta_0}{2},$$
for $\theta_0\leq\theta\leq \pi,\; z\in S_{R}$ and $R>R_0$,
and this yields (\ref{4s}).

It follows from (\ref{4s}) that
\begin{eqnarray*}\ds
&{}&\left|\int_{\theta_0\leq|\theta|\leq\pi}
\exp\left(kQ_k(F_k(z)+r_k(z)e^{i\theta})-ik\theta\right)d\theta\right|\\ \\
&\leq&\int_{\theta_0\leq|\theta|\leq\pi}
\exp\left\{\Re\left(kQ_k(F_k(z)+r_k(z)e^{i\theta})\right)\right\}d\theta\\ \\
&\leq& 2(\pi-\theta_0)\exp\left\{ k\left(\Rea Q_k(z)-\frac12(1-\cos\theta_0)\right)\right\}
\\ \\
&=&o\left(\frac{e^{kQ_k(z)}}{\sqrt{k}}\sqrt{\frac{\pi}{C(z)}}\right).
\end{eqnarray*}
Combining this with (\ref{4e1}) we see that
$$\ds
\int_{|\theta|\leq\pi}\exp\left(kQ_k(F_k(z)+r_k(z)e^{i\theta})
-ik\theta\right)d\theta\sim\frac{e^{kQ_k(z)}}{\sqrt{k}}\sqrt{\frac{\pi}{C(z)}}$$
as $k\to\infty$, uniformly on compact subsets in $S_R$, where $R>R_0$.
Together with (\ref{4a1}) this gives
$$\ds
f^{(k)}(a_kF_k(z))\sim\frac{k!q_k^k(z)e^{b_k}e^{kQ_k(z)}}{
2a_k^k\sqrt{\pi k C(z)}}.$$
Taking logarithms, dividing by $k$ and passing to the limit
as $k\to\infty$ we obtain using (\ref{defuk})
$$u(F(z))=\Rea Q(z)+\log|q(z)|,\quad z\in S_{R},$$
for $R>R_0$.
This is the same as (\ref{identity}). One consequence of (\ref{identity}) is that
$u\not\equiv-\infty$.

\vspace{.2in}

\noindent
{\bf 5. Analytic continuation of $u$ and conclusion of the proof}
\nopagebreak
\vspace{.2in}

As pointed out at the end of Section~3, in order to prove 
Theorem 1 it suffices to show that the function $u$ obtained
is not harmonic.
\vspace{.1in}

\noindent
{\bf Lemma 6.} {\em Let
\begin{equation}
\label{26}
q(z)=-z^{2p}\psi(z),
\end{equation}
where $p\geq 1$ and $\psi$ is an analytic function
mapping $H$ to $\overline{H}\setminus \{0\}$.
Define
$$Q(z):=\int_i^zq(\zeta)d\zeta$$
and
\begin{equation}
\label{29}
F(z):=z-1/q(z).
\end{equation}
Let $u$ be a subharmonic function in $H$ satisfying $(\ref{identity})$,
in a region $S_R$ defined in
$(\ref{defSR})$.
Then $u$ is not harmonic in $H$.}
\vspace{.1in}

{\em Proof.}
Suppose that $u$ is harmonic.
Then there exists a holomorphic function $h$ in $H$ such that
$u=\Rea h$, and 
\begin{equation}
\label{dur}
h(F(z))=Q(z)+\log q(z),\quad z\in S_{R}.
\end{equation}
Differentiating (\ref{dur})
and using $q=Q'$ and (\ref{29})
we obtain
$$h'(F(z))=q(z)=\frac{1}{z-F(z)}.$$
{}From (\ref{26}), (\ref{29}) and (\ref{21}), it follows that there is
a branch $G$ of the inverse $F^{-1}(w)$ which is defined 
in a Stolz angle $S$ 
and satisfies
\begin{equation}\label{kkk}
G(w)\sim w\quad\mbox{as}\quad w\to\infty,\; w\in S.
\end{equation}
In particular,
$G(w)\in H$ for $w\in S$ and $|w|$ large enough.
We have
\begin{equation}
\label{fa}
h'(w)=q(G(w))=\frac{1}{G(w)-w},
\end{equation}
for $w\in S$ and $|w|$ large enough.
Since $h$ is holomorphic in
$H$ we see that $G$ has a meromorphic continuation to $H$.
Using (\ref{26}),
the second equation in (\ref{fa}) can  be rewritten as
\begin{equation}
\label{5b}
\psi(G(w))
=\frac{1}{G^{2p}(w)(w-G(w))}.
\end{equation}
We will derive from (\ref{5b}) that $G$ maps $H$ into itself.
To show this, we establish first that $G$ never takes a real
value in $H$. In view of (\ref{kkk}),
there exists a point $w_0\in H$
such that $G(w_0)\in H$. If $G$ takes a real value in $H$, then there
exists a curve $\phi:[0,1]\to H$ beginning at $w_0$ and ending at some
point $w_1\in H$, such that $G(\phi(t))\in H$ for $0\leq t<1$
but $G(w_1)=G(\phi(1))\in \R$.
We may assume that $G(w_1)\neq 0$; this can be achieved by
a small perturbation of the curve $\phi$ and the point $w_1$.
Using (\ref{5b}) we obtain an analytic continuation
of $\psi$ to the real point $G(w_1)$ along the curve $G(\phi)$.
We have 
$$\lim_{t\to 1}\Ima\psi(G(\phi(t)))\geq 0,$$
because the imaginary part of $\psi$ is non-negative in $H$.

It follows that for $w\to w_1$, the right hand side of (\ref{5b}) 
has negative imaginary part, while the left hand side has
non-negative imaginary part, which is a contradiction.
Thus we have proved that $G$ never takes real values in~$H$. 

Since $G(w_0)\in H$ we see that $G$ maps $H$ into itself.
Then Lemma~3
and (\ref{kkk}) imply that
$\Ima (G(w)-w)>0\quad\mbox{for}\quad w\in H.$
Combining this with the second equation
of (\ref{fa}) we obtain $\Ima q(G(w))<0$ for $w\in H$.
Using (\ref{kkk}) we find that in particular $\Ima q(e^{i\pi/(2p)}y)<0$
for large $y>0$.
On the other hand, we have
$$\Ima q(e^{i\pi/(2p)}y)=y^{2p}\Ima \psi(e^{i\pi/(2p)}y)\geq 0$$
by (\ref{26}).
This contradiction proves Lemma~6, and it also completes the proof of Theorem~1.
\vspace{.1in}

{\em Proof of Theorem 2.} 
It remains to show that there exist 
positive constants $\gamma$ and $\alpha$ depending only on $p$,
such that the total Riesz measure of $u$
in the region  $\{ z: \Ima z> \gamma |z|\}$ is 
bounded below by $\alpha/2$,
with $u$ and $p$ as in Sections 3 and 4.

First we note that 
\begin{equation}
\label{boundu}
u(z)\leq\max_{|\zeta|=t}\Rea Q(z+\zeta)-\log t, \quad t=\frac12 \Ima z,
\end{equation}
for $z\in H$ by (\ref{bounduk}).
Using (\ref{defQ}), (\ref{4a}), (\ref{i}) and (\ref{21}) we see that
for every compact subset $K$ of $H$, the right hand side
of (\ref{boundu}) is bounded from above by a constant that depends only on
$p$ and $K$.

Thus for fixed $p$, the functions $u$, $Q$ and $q$
under consideration belong to normal
families. 
Arguing by contradiction, we assume
that 
$\alpha$ and $\gamma$ as above do not exist.
Then there exists 
a sequence $(u_k)$ of subharmonic functions in $H$, satisfying equations
of the form
$$u_k(F_k)=\Rea Q_k+\log|q_k|$$
similar to (\ref{identity}), such that the Riesz measure of $u_k$ in
$\{ z: \Ima z> |z|/k\}$ tends to $0$ as $k\to\infty$.
(The functions $u_k, F_k, q_k, Q_k$ are not the functions introduced
in Sections~3 and~4.)
It is important that all these equations hold in the same region $S_R$.
Using normality we can take convergent subsequences,
and we obtain a limit equation of the same form, satisfied
by a function $u$ harmonic in $H$. 
But this contradicts Lemma~6, and thus the proof of Theorem~2 is completed.
\vspace{.15in}

\noindent
{\bf Remark.}
In \cite{SS} and in most subsequent work on the subject the
auxiliary function $F_0(z)=z-f(z)/f'(z)$ plays an important role.
Note that our function $F$ of (\ref{29}) is of a similar nature, except that
$f'/f$ is replaced by $q$ which is obtained from $f'/f$ by rescaling.

{\em W. B.: Mathematisches Seminar,
Christian-Albrechts-Universit\"at zu Kiel,
Lude\-wig-Meyn-Str.~4,
D-24098 Kiel,
Germany

bergweiler@math.uni-kiel.de

A. E.: Department of Mathematics,
Purdue University, West Lafayette, IN 47907, USA

eremenko@math.purdue.edu}

\end{document}